\documentclass[12pt]{amsart}
\usepackage{graphicx}
\usepackage{amssymb}
\usepackage{amsmath}
\usepackage{amscd}
\usepackage{upgreek}
\usepackage{mathrsfs}
\usepackage{stmaryrd}
\usepackage{longtable}
\usepackage{txfonts}
\usepackage[T1]{fontenc}
\usepackage{eucal}
\usepackage{titletoc}
\usepackage{wrapfig}
\usepackage{float}
\usepackage{xypic}
\usepackage{dsfont}
\usepackage{xcolor}
\usepackage{color}
\usepackage[colorlinks,linkcolor=blue,anchorcolor=blue,urlcolor=blue,
citecolor=blue]{hyperref}
\usepackage{hyperref}
\usepackage[hmarginratio=1:1,top=32mm,columnsep=20pt]{geometry}
\DeclareMathOperator{\Hom}{Hom}	
\allowdisplaybreaks	
\usepackage{mathptmx}
\DeclareMathAlphabet{\mathcal}{OMS}{cmsy}{m}{n}
\DeclareSymbolFont{largesymbols}{OMX}{cmex}{m}{n}

\newtheorem{theorem}{Theorem}[section]
\newtheorem{prop}[theorem]{Proposition}
\newtheorem{defn}[theorem]{Definition}
\newtheorem{lem}[theorem]{Lemma}
\newtheorem{coro}[theorem]{Corollary}

\newtheorem{thm}[theorem]{Theorem}

\newtheorem{rem}[theorem]{\rm \textsc{Remark}}
\newtheorem{exam}[theorem]{Example}
\newtheorem{exe}[theorem]{\rm \textsc{Exercise}}

\newcommand{\C}{\mathds{C}}
\newcommand{\N}{\mathbb{N}}
\newcommand{\Q}{\mathbb{Q}}
\newcommand{\Z}{\mathds{Z}}

\newcommand{\ra}{\longrightarrow}

\newcommand{\tr}{{\rm tr}}

\newcommand{\vp}{\varphi}

\begin{document}
\setlength{\oddsidemargin}{0cm}
\setlength{\evensidemargin}{0cm}

%%%%%%%%%%%%%%%%%%%%%%%%%%%%%%%%%%% TITLE, AUTHOR AND ADDRESS
\title{\scshape Bogomolov Multipliers of Some Groups of order $p^{6}$}
\author{\scshape Yin Chen}

\address{Department of Mathematics, Nantong University, Nantong 226019, P.R. China}

\email{ychen@ntu.edu.cn}

\author{\scshape Rui Ma}
\address{School of Mathematics and Statistics, Northeast Normal University, Changchun 130024, P.R. China}
\email{mar770@nenu.edu.cn}
\date{\today}
\def\shorttitle{Bogomolov Multipliers of Some Groups of order $p^{6}$}

%%%%%%%%%%%%%%%%%%%%%%%%%%%%%%%%%%%ABSTRACT, MSC AND KEYWORDS
\begin{abstract}
Let $G$ be a finite group, $V$ a faithful finite-dimensional representation of $G$ over the complex field $\C$ and $\C(V)^{G}$ be the corresponding invariant field.  The Bogomolov multiplier $B_{0}(G)$ of $G$ is canonically isomorphic to the unramified cohomological group $H_{\textrm{nr}}^{2}(\C(V)^{G},\Q/\Z)$, which has been used by Saltman (1984) and Bogomolov (1988) to provide counter-examples to the rationality problem of  $\C(V)^{P}$ for finite $p$-groups $P$ over $\C$. In this paper, we investigate the vanishing  property of $B_{0}(P)$, where $P$ denotes a $p$-group of order $p^{6}$ for $p\geqslant3$.
\end{abstract}


\maketitle
\baselineskip=18pt

%%%%%%%%%%%%%%%%%%%%%%%%%%%%%%%%%%%%%CONTENTS
%\textcolor{blue}{\tableofcontents{}}

\dottedcontents{section}[1.16cm]{}{1.8em}{5pt}
\dottedcontents{subsection}[2.00cm]{}{2.7em}{5pt}
\dottedcontents{subsubsection}[2.86cm]{}{3.4em}{5pt}

%%%%%%%%%%%%%%%%%%%%%%%%%%%%%%%%%%%%%%%%%%%SECTION 1
\section{\scshape Introduction}
\setcounter{equation}{0}
\renewcommand{\theequation}
{1.\arabic{equation}}
\setcounter{theorem}{0}
\renewcommand{\thetheorem}
{1.\arabic{theorem}}

Let $G$ be a finite group and $H^{2}(G,\Q/\Z)$ be the Schur multiplier of $G$. For a subgroup $A\leqslant G$, we use
$\textrm{res}_{A}^{G}$ to denote the usual restriction map from $H^{2}(G,\Q/\Z)$  to $H^{2}(A,\Q/\Z)$. The \textbf{Bogomolov multiplier} $B_{0}(G)$ of $G$ is defined as the intersection of all kernels of $\textrm{res}_{A}^{G}$ where $A$ runs over all bicyclic subgroups of $A$. In 1988, Bogomolov \cite{Bog1988} proved that $B_{0}(G)$ is canonically isomorphic to the unramified cohomological group $H_{\textrm{nr}}^{2}(\C(V)^{G},\Q/\Z)$, which is important to the birational geometry of the quotient variety $V/G$, especially to answer the (retract) rationality problem of the invariant field $\C(V)^{G}$; see Artin-Mumford \cite{AM1972} and Saltman \cite{Sal1984}. In particular, Bogomolov \cite{Bog1988}  used the non-vanishing property of $B_{0}(P)$ as an obstruction to find a $p$-group $P$ of order $p^{6}$ such that the classical Noether's problem over $\C$ has a negative answer, i.e., the invariant field $\C(P):=\C(x_{g}\mid g\in P)^{P}$ is not rational over $\C$. This reduces the order of the finite $p$-group in the first counter-example of Saltman \cite{Sal1984} to the classical Noether's problem over $\C$. Since then, it is natural to classify nonabelian
finite $p$-groups $P$ of order $\leqslant p^{6}$ with $B_{0}(P)=0$; see Bogomolov \cite[Remark 1, page 479]{Bog1988}.
Note that the classical Noether's problem for finite abelian groups has been answered completely; see for example Swan \cite{Swa1983} for a survey.

A result due to Chu-Kang \cite{CK2001} states that for a prime $p$ and any $p$-group $P$ of order $\leqslant p^{4}$, the invariant field $\C(V)^{P}$ is rational over $\C$, where $V$ denotes a faithful finite-dimensional representation of $G$ over $\C$.
This result particularly implies that $B_{0}(P)=0$ for any $p$-group $P$ of order $\leqslant p^{4}$.
In \cite{Bog1988}, Bogomolov also claimed that $B_{0}(P)=0$ for any $p$-group $P$ of order $p^{5}$.
Chu-Hu-Kang-Prokhorov \cite{CHKP2008} affirmatively answered Noether's problem for the groups of order $32$, thus confirming Bogomolov's claim for the case $p=2$. However, Moravec \cite{Mor2012b} provided three groups of
order $3^{5}$ with nonzero Bogomolov multipliers to demonstrate that Bogomolov's claim was not correct for any prime $p>2$.
Furthermore, Moravec \cite{Mor2014} proved that if two finite groups $G$ and $H$ are isoclinic, then $B_{0}(G)$ is isomorphic to
$B_{0}(H)$. This result leads one to define the Bogomolov multiplier $B_{0}(\Phi):=B_{0}(G)$, where $\Phi$ is an isoclinic family containing a group $G$. Recall that for a prime $p\geqslant 3$, James \cite{Jam1980} in 1980 has already classified all nonabelian  groups of order $p^{5}$ and $p^{6}$ into 9 isoclinic families: $\Phi_{2},\dots, \Phi_{10}$ and into 42  isoclinic families: $\Phi_{2},\dots, \Phi_{43}$ respectively.
Moravec \cite{Mor2012} and Hoshi-Kang-Kunyavskii \cite{HKK2013} used different methods to prove that if $P$ is a nonabelian group of order $p^{5} (p\geqslant 3)$, then $B_{0}(P)=0$ if and only if $P\notin \Phi_{10}$. In 2010,
Chu-Hu-Kang-Kunyavskii \cite{CHKK2010} classified all nonabelian groups $P$ of order $2^{6}$ with $B_{0}(P)=0$. Recently, Michailov \cite{Mic2014} and \cite{Mic2016} also studied the Bomogolov multiplier of some special kinds of groups of order $p^6$.

The purpose of this article is to study the vanishing property of $B_{0}(P)$ for a nonabelian group $P$ of order $p^{6} (p\geqslant3)$. Our main result can be summarized as follows.

\begin{thm}\label{main-thm}
Let $p\geqslant3$ be a prime and $\{\Phi_{k}\mid 2\leqslant k\leqslant 43\}$ be the set of isoclinic families of groups of order $p^{6}$. Then $B_{0}(\Phi_{k})$ is not zero for $k\in \Delta:=\{10,18,20,21,36,38,39\}$ and $B_{0}(\Phi_{i})=0$ for $k\in \{2,3,\dots,43\}\setminus (\Delta \cup\{15,28,29\}).$
\end{thm}

\begin{rem}{\rm
With the help of a computer calculation, it was conjectured that  $B_{0}(\Phi_{15})=B_{0}(\Phi_{28})=B_{0}(\Phi_{29})=0$.
However, we currently are not able to confirm this conjecture.
}\end{rem}

We notice that there are a number of articles addressing the vanishing property of the Bogomolov multipliers for non $p$-groups; see, for example, finite simple groups (Bogomolov-Maciel-Petrov \cite{BMP2004} and
Kunyavskii \cite{Kun2010}), unitriangular groups (Michailov \cite{Mic2015}), and rigid finite groups
(Kang-Kunyavskii \cite{KK2014}, Rai-Yadav \cite{RY2015}, and Rai \cite{Rai2019}).

In Section 2, we recall a group-theoretic description due to Moravec \cite{Mor2012b} for $B_{0}(G)$ and collect some
well-known identities and properties on nilpotent groups, commutator subgroups and nonabelian exterior squares, which will be used repeatedly in our proofs.
Section 3 is devoted to proving, case by case, that  $B_{0}(\Phi_{i})=0$ for $k\in \{2,3,\dots,43\}\setminus (\Delta \cup\{15,28,29\}).$
In Section 4, we prove that $B_{0}(\Phi_{k})$ is not zero for $k\in \Delta$.
Throughout this article, for a group $G$ and $x,y\in G$, we define $x^{y}:=y^{-1}xy$ and write
$x^{-1}x^{y}=x^{-1}y^{-1}xy$ for the commutator $[x,y]$. We also define $[x_{1},\dots,x_{n}]:=[\cdots[x_{1},\dots,x_{n-1}],x_{n}]$ for $x_{1},\dots,x_{n}\in G$ and $n\geqslant 2$. In particular, we use $[x,_{n}y]$ to denote $[x,y,\dots,y]$ with $n$ copies of $y$.

\section{\scshape Preliminaries}
\setcounter{equation}{0}
\renewcommand{\theequation}
{2.\arabic{equation}}
\setcounter{theorem}{0}
\renewcommand{\thetheorem}
{2.\arabic{theorem}}

Let $G$ be a group and $\vp$ be an automorphism of $G$. We denote by $G^{\vp}:=\{x^{\vp}\mid x\in G\}$ the
isomorphic copy of $G$ via $x\mapsto x^{\vp}$. Define $\tau(G)$ to be the group generated by
$G$ and $G^{\vp}$, subject to the commutator relations:
$$[x,y^{\vp}]^{z}=[x^{z},(y^{z})^{\vp}]=[x,y^{\vp}]^{z^{\vp}}\textrm{ and }[x,x^{\vp}]=1$$
for all $x,y,z\in G$. The groups $G$ and $G^{\vp}$ can be regarded to subgroups of $\tau(G)$.
Let $[G,G^{\vp}]:=\langle [x,y^{\vp}]\mid x,y\in G\rangle$ be the commutator subgroup in $\tau(G)$.
Note that $[G,G^{\vp}]$ is isomorphic to $G\wedge G$, the nonabelian exterior square of $G$; see for example Blyth-Morse \cite[Proposition 16]{BM2009}.

Consider the commutator subgroup $[G,G]$ of $G$.
There exists a natural group homomorphism $\rho: [G,G^{\vp}]\ra [G,G]$
defined by $[x,y^{\vp}]\mapsto [x,y]$ for all $x,y\in G$. Let $M(G)$ be the kernel of $\rho$, i.e.,
$$M(G):=\left\{\prod_{\textrm{finite}}[x_{i},y_{i}^{\vp}]^{\epsilon_{i}}\in [G,G^{\vp}]~\Big |~ \epsilon_{i}=\pm 1, \prod_{\textrm{finite}}[x_{i},y_{i}]^{\epsilon_{i}}=1\right\}$$
and define
$$M_{0}(G):=\left\{\prod_{\textrm{finite}}[x_{i},y_{i}^{\vp}]^{\epsilon_{i}}\in [G,G^{\vp}]~\Big |~ \epsilon_{i}=\pm 1, [x_{i},y_{i}]=1,\textrm{ for each }i\right\}.$$
Clearly, $M_{0}(G)$ is a subgroup of $M(G)$. Furthermore, $B_{0}(G)$ is isomorphic to the quotient group $M(G)/M_{0}(G)$; see
Moravec \cite[Section 2]{Mor2012}.

Here we collect some useful properties of $\tau(G)$ and $[G,G^{\vp}]$ as follows.

\begin{prop}\label{prop2.1}
Let $G$ be a group, $\vp$ be an automorphism of $G$ and $x,y,z,w\in G$ be arbitrary elements. Then
\begin{enumerate}
  \item $[x,yz]=[x,z][x,y]^{z}=[x,z][x,y][x,y,z]$ and $[xy,z]=[x,z]^{y}[y,z]=[x,z][x,z,y][y,z]$.
  \item $[x,y^{\vp}]=[x^{\vp},y]$.
  \item $[x,y,z^{\vp}]=[x,y^{\vp},z]=[x^{\vp},y,z]=[x,y^{\vp},z^{\vp}]=[x^{\vp},y^{\vp},z]=[x^{\vp},y,z^{\vp}]$.
  \item $[[x,y^{\vp}],[z,w^{\vp}]]=[[x,y],[z,w]^{\vp}]$.
  \item If $[x,y]=1$, $[x^{n},y^{\vp}]=[x,y^{\vp}]^{n}=[x,(y^{\vp})^{n}]$ for all $n\in \Z$.
\end{enumerate}
\end{prop}

\begin{proof}
The equalities in $(1)$ can be verified through direct computations. For the remaining statements, we note that
the nonabelian exterior square $G\wedge G$  of  a group $G$ is a quotient of the nonabelian tensor product $G\otimes G$, thus
the equalities in $(2)$--$(5)$ follow from Blyth-Morse \cite[Lemmas 9, 10, and 11]{BM2009}.
\end{proof}

Given a group $G$, we define $G^{1}:=G$ and $G^{n}:=[G^{n-1},G]$ for $n\geqslant 2.$ Recall that $G$ is \textbf{nilpotent} if
there exists a $c\in\N$ such that $G^{c+1}=\{1\}$; the least such $c$ is called the \textbf{class} of the nilpotent group $G$.

\begin{prop}\label{prop2.2}
Let $G$ be a group and $\vp$ be an automorphism of $G$.
\begin{enumerate}
  \item If $G$ is nilpotent of class $c$, then $\tau(G)$ is nilpotent of class $\leqslant c+1$.
  \item If $[G,G]$ is nilpotent of class $c$, then $[G,G^{\vp}]$ is nilpotent of class $c$ or $c+1$.
  \item If $G$ is nilpotent of class $\leqslant 2$, then $[G,G^{\vp}]$ is abelian.
  \item If $G$ is nilpotent group of class $\leqslant 3$, then
$$[x,y^{n}]=[x,y]^{n}[x,_{2}y]^{n\choose 2}[x,_{3}y]^{n\choose 3}$$
for all $x,y\in \tau(G)$ and all $n\in \N^{+}$.
\item If $H$ is nilpotent group of class $\leqslant 5$, then
$$[x^{n},y]=[x,y]^{n}[x,y,x]^{n\choose 2}[x,y,_{2}x]^{n\choose 3}[x,y,_{3}x]^{n\choose 4}[x,y,x,[x,y]]^{\delta(n)}$$
for all $x,y\in H$ and all $n\in \N^{+}$, where $\delta(n):=\frac{n(n-1)(2n-1)}{6}$. In particular, if $[x,y]$ belongs to the center of $H$, then $[x^{n},y]=[x,y]^{n}$.
\end{enumerate}
\end{prop}

\begin{proof}
These statements follow from Moravec \cite[Lemmas 2.1, 3.1, and 3.7]{Mor2012}.
\end{proof}

Recall that a group $G$ is called \textbf{polycyclic} if there exists a subnormal series
$G=G_{1}\vartriangleright G_{2}\vartriangleright \dots\vartriangleright G_{n+1}=\{1\}$ such that each factor $G_{k}/G_{k+1}$
is cyclic of order $r_{k}$. Clearly, a polycyclic group is solvable.
A sequence $x_{1},\dots,x_{n}$ of elements of a finite polycyclic group $G$  is called a \textbf{polycyclic generating sequence} of $G$ if $G_{k}$ can be generated by $G_{k+1}$ and $x_{k}$ for $k=1,\dots,n$. The number $r_{k}$ is called the \textbf{relative order} of $x_{k}$. An element $x$ in a polycyclic generating sequence is said to be \textbf{absolute} if
the relative order of $x$ and the order of $x$ are equal.

\begin{lem}[\cite{BM2009}, Proposition 20]\label{lem2.3}
Let $G$ be a finite polycyclic group with a polycyclic generating sequence $x_{1},\dots,x_{n}$. Then
$[G,G^{\varphi}]$, as a subgroup of $\tau(G)$, is generated by $\{[x_{i},x_{j}^{\vp}]\mid 1\leqslant i<j\leqslant n\}$.
\end{lem}

\begin{lem}[\cite{Mor2012}, Proposition 3.2]\label{lem2.4}
Let $p\geqslant 3$ be a prime and $G$ be a finite $p$-group of nilpotency class $\leqslant 3$.
Let $x_{1},\dots,x_{n}$ be a polycyclic generating sequence of $G$. If all nontrivial commutators $[x_{i},x_{j}](i<j)$
are different absolute elements of $\{x_{1},\dots,x_{n}\}$, then $B_{0}(G)=0$.
\end{lem}

%\begin{rem}{\rm As the vanishing property of the Bogomolov multipliers of groups of order $3^{5}$ can be obtained by computer calculations, Moravic assumed throughout the whole article \cite{Mor2012} that the prime $p\geqslant 5$. With a careful reading, we see that  \cite[Proposition 3.2]{Mor2012} follows for all prime $p\geqslant 3$.
%}\end{rem}

%\begin{lem}\label{lem2.5}
%Let $p\geqslant 5$ be a prime and $G$ be a finite $p$-group of nilpotency class $\leqslant 4$. Let $x_{1},\dots,x_{n}$ be a polycyclic generating sequence of $G$. Suppose all nontrivial commutators $[x_{i},x_{j}](i<j)$ are different absolute elements of $\{x_{1},\dots,x_{n}\}$. If $[G,G^{\varphi}]$ has nilpotency class $\leqslant 2$, then $B_{0}(G)=0$.
%\end{lem}

%\begin{proof} $[x,y]^{p}=1$. $[x^{p},y']=[x,y']^{p}[x,y',x]^{p\choose 2}[x,y',_{2}x]^{p\choose 3}[x,y',_{3}x]^{p\choose 4}[x,y',x,[x,y']]^{\delta(p)}$
%\begin{eqnarray*}
%[y',x]^{p}& = & [y',x^{p}] [x,y',x]^{p\choose 2}[x,y',_{2}x]^{p\choose 3}[x,y',_{3}x]^{p\choose 4}[x,y',x,[x,y']]^{\delta(p)}.
%\end{eqnarray*}

%\begin{eqnarray*}
%[x',[x,y]]^{n}& = &  [x',[x,y]^{p}] [[x,y],x',[x,y]]^{p\choose 2}[[x,y],x',_{2}[x,y]]^{p\choose 3}[[x,y],x',_{3}[x,y]]^{p\choose 4}[[x,y],x',[x,y],[[x,y],x']]^{\delta(p)}.
%\end{eqnarray*}

%\end{proof}

\begin{lem}[\cite{Kan2014}, Theorem 1.4]\label{lem2.6}
 Let $G$ and $H$ be two finite groups. Then $B_{0}(G\times H)$ is isomorphic to $B_{0}(G)\times B_{0}(H)$.
\end{lem}

\section{\scshape Vanishing Bogomolov Multipliers}
\setcounter{equation}{0}
\renewcommand{\theequation}
{3.\arabic{equation}}
\setcounter{theorem}{0}
\renewcommand{\thetheorem}
{3.\arabic{theorem}}

Throughout this and next sections, we use $\nu$ to denote the smallest positive integer which is non-quadratic residue modulo $p$ and use $\mu$ to denote
the smallest positive integer which is a primitive root modulo $p$; we define $\alpha_{i+1}^{(p)}:=\alpha_{i+1}^{p}\alpha_{i+2}^{p\choose 2}\cdots \alpha_{i+k}^{p\choose k}\cdots \alpha_{i+p}^{p\choose p}$; we also follow James's notations in the list of isoclinic families of groups of orders $p^{5}$ and $p^{6}$ for $p\geqslant3$; see \cite{Jam1980}.

\begin{prop}\label{prop3.1}
$B_{0}(\Phi_{k})=0$, for $k\in \{2,3,\dots,17\}\setminus\{10,15\}$ and $B_{0}(\Phi_{10})\neq 0$.
\end{prop}

\begin{proof}
It follows from Lemma \ref{lem2.6} that $B_{0}(\Phi_{2})=B_{0}(\Phi_{2}(411)_{a})\cong B_{0}(\Phi_{2}(41))\times B_{0}((1))$, where
$\Phi_{2}(41)$ belongs to the second isoclinic family of groups of order $p^{5}$ and $(1)$ denotes the cyclic group of order $p$.
By Hoshi-Kang-Kunyavskii \cite[Theorem 1.12]{HKK2013} we see that $B_{0}(\Phi_{2}(41))=0$. Further, it follows from Chu-Kang \cite[Theorem 1.6]{CK2001} that $B_{0}((1))=0$. This proves $B_{0}(\Phi_{2})=0$. Applying the same method, we see that
$B_{0}(\Phi_{3})=B_{0}(\Phi_{4})=\dots=B_{0}(\Phi_{9})=B_{0}(\Phi_{12})=0$ and $B_{0}(\Phi_{10})\neq 0$.

Note that $\Phi_{11}(1^{6})\in \Phi_{11}$ has polycyclic presentation: $\{\alpha_{1},\beta_{1},\alpha_{2},\beta_{2},\alpha_{3},\beta_{3}\}$ is the generating set, subject to all nontrivial commutator relations $[\alpha_{1},\alpha_{2}]=\beta_{3},[\alpha_{2},\alpha_{3}]=\beta_{1}, [\alpha_{3},\alpha_{1}]=\beta_{2}$ and the orders of all generators are $p$. As $\Phi_{11}(1^{6})$ has nilpotency class 2, it follows Lemma \ref{lem2.4} that $B_{0}(\Phi_{11})=B_{0}(\Phi_{11}(1^{6}))=0$. Similarly, we have
$B_{0}(\Phi_{14})=B_{0}(\Phi_{16})=B_{0}(\Phi_{17})=0.$

To show $B_{0}(\Phi_{13})=0$, we take $G=\Phi_{13}(1^{6})$ as a representative in $\Phi_{13}$.
Comparing with the polycyclic presentation of $G$ given by James \cite[page 626]{Jam1980} and replacing $\alpha_{3}$ by $\alpha_{3}^{p-1}$, we obtain a new polycyclic presentation for $G$:
$$\langle \alpha_{1},\dots,\alpha_{4},\beta_{1},\beta_{2}\mid [\alpha_{1},\alpha_{2}]=\beta_{1},
[\alpha_{3},\alpha_{1}]=\beta_{2}=[\alpha_{2},\alpha_{4}],\alpha_{1}^{p}=\alpha_{2}^{p}=\alpha_{3}^{p}=\alpha_{4}^{p}=\beta_{1}^{p}=\beta_{2}^{p}=1\rangle$$
where the trivial commutator relations have been deleted. As $G$ has nilpotency class 2, Proposition \ref{prop2.2} (3) implies that $[G,G^{\vp}]$ is abelian.  It follows from Lemma \ref{lem2.3} that the group $[G,G^{\vp}]$ is generated by $[\alpha_{1},\alpha_{2}^{\vp}],[\alpha_{1},\alpha_{3}^{\vp}]$ and $[\alpha_{2},\alpha_{4}^{\vp}]$ modulo $M_{0}(G)$, which implies that each $w\in [G,G^{\vp}]$ can be written as
$$w=[\alpha_{1},\alpha_{2}^{\vp}]^{r}[\alpha_{1},\alpha_{3}^{\vp}]^{s}[\alpha_{2},\alpha_{4}^{\vp}]^{t}\cdot w_{0}$$
for some $w_{0}\in M_{0}(G)$. Now suppose $w\in M(G)$ is an arbitrary element. To show $B_{0}(G)=0$, it suffices to show $w\in M_{0}(G)$. As $1=\rho(w)=\beta_{1}^{r}\beta_{2}^{t-s}$ and $\beta_{1},\beta_{2}$ are elements in the polycyclic generating sequence, we see that $\beta_{1}^{r}=\beta_{2}^{t-s}=1$. Note that $\beta_1$ and $\beta_2$ have order $p$. Hence, $r$ and $t-s$ both are divisible by $p$.
We \textbf{claim} that $[\alpha_{1},\alpha_{2}^{\vp}]^{p}=1$. By Proposition \ref{prop2.2} (4), we have
$$1=[\alpha_{1}^{\vp},1]=[\alpha_{1}^{\vp},\alpha_{2}^{p}]=[\alpha_{1}^{\vp},\alpha_{2}]^{p}[\alpha_{1}^{\vp},_{2}\alpha_{2}]^{p\choose 2}[\alpha_{1}^{\vp},_{3}\alpha_{2}]^{p\choose 3}.$$
Proposition \ref{prop2.2}  (1) implies that $\tau(G)$ has nilpotency class at most 3.
Thus $[\alpha_{1}^{\vp},_{3}\alpha_{2}]=1$. Further, by Proposition \ref{prop2.1}  (3), we see that
 $[\alpha_{1}^{\vp},_{2}\alpha_{2}]^{p\choose 2}=[\alpha_{1}^{\vp},\alpha_{2},\alpha_{2}]^{p\choose 2}=
 [\alpha_{1},\alpha_{2},\alpha_{2}^{\vp}]^{p\choose 2}= [\beta_{1},\alpha_{2}^{\vp}]^{p\choose 2}$. Since $[\beta_{1},\alpha_{2}]=1$, it follows from Proposition \ref{prop2.1}  (5) that $ [\beta_{1},\alpha_{2}^{\vp}]^{p\choose 2}= [\beta_{1}^{p\choose 2},\alpha_{2}^{\vp}]=[1,\alpha_{2}^{\vp}]=1.$
Hence, by Proposition \ref{prop2.1} (2), we have  $[\alpha_{1},\alpha_{2}^{\vp}]^{p}=[\alpha_{1}^{\vp},\alpha_{2}]^{p}=1$ and the claim follows.
A similar argument shows that $[\alpha_{1},\alpha_{3}^{\vp}]^{p}=1=[\alpha_{2},\alpha_{4}^{\vp}]^{p}$. Thus $w$ can be written as
$$w=([\alpha_{1},\alpha_{3}^{\vp}][\alpha_{2},\alpha_{4}^{\vp}])^{s}w_{0}.$$
Now it suffices to show $[\alpha_{1},\alpha_{3}^{\vp}][\alpha_{2},\alpha_{4}^{\vp}]\in M_{0}(G)$.
Note that $G$ has nilpotency class 2 and by Proposition \ref{prop2.1}  (1), we have
$[\alpha_1\alpha_2,\alpha_3\alpha_4]=[\alpha_1,\alpha_4][\alpha_2,\alpha_4][\alpha_1,\alpha_3]
[\alpha_2,\alpha_3]=1.$ This means that $[\alpha_1\alpha_2,\alpha_3^\varphi\alpha_4^\varphi]\in M_0(G)$. Expanding this element by repeatedly applying the identities in Proposition \ref{prop2.1}, we  obtain
\begin{eqnarray*}
[\alpha_1\alpha_2,\alpha_3^\varphi\alpha_4^\varphi]
&=&[\alpha_1\alpha_2,\alpha_4^\varphi]\cdot[\alpha_1\alpha_2,\alpha_3^\varphi]\cdot
[\alpha_1\alpha_2,\alpha_3^\varphi,\alpha_4^\varphi]    \\
&=& [\alpha_1,\alpha_4^\varphi] [\alpha_1,\alpha_4^\varphi,\alpha_2][\alpha_2,\alpha_4^\varphi]\cdot
[\alpha_1,\alpha_3^\varphi] [\alpha_1,\alpha_3^\varphi,\alpha_2][\alpha_2,\alpha_3^\varphi]\cdot\\
&&
[[\alpha_1,\alpha_3] [\alpha_1,\alpha_3,\alpha_2][\alpha_2,\alpha_3],\alpha_4^\varphi]\\
&=& [\alpha_1,\alpha_4^\varphi][\alpha_2,\alpha_4^\varphi]\cdot
[\alpha_1,\alpha_3^\varphi] [\beta_2^{-1},\alpha_2^\varphi][\alpha_2,\alpha_3^\varphi]\cdot [\beta_2^{-1},\alpha_4^\varphi].
\end{eqnarray*}
Since $[\alpha_1,\alpha_4^\varphi],[\beta_2^{-1},\alpha_2^\varphi],[\alpha_2,\alpha_3^\varphi], [\beta_2^{-1},\alpha_4^\varphi]$ and $[\alpha_1\alpha_2,\alpha_3^\varphi\alpha_4^\varphi]$ are elements in $M_0(G)$, it follows that
$[\alpha_{1},\alpha_{3}^{\vp}][\alpha_{2},\alpha_{4}^{\vp}]\in M_{0}(G)$. Hence, $B_{0}(\Phi_{13})=0$ and the proof is completed.
\end{proof}

\begin{rem}{\rm
Noether's problem for groups of order $p^6$ that have an abelian normal subgroup such that the quotient group
is cyclic has been settled by Michailov \cite[Theorem 1.9]{Mic2014}. Using this result we also can obtain 
that $B_{0}(\Phi_{14})=0$.
}\end{rem}

\begin{prop}
$B_{0}(\Phi_{19})=0$.
\end{prop}

\begin{proof}
Let $G=\Phi_{19}(1^6)$. Then it is generated by $\{\alpha,\alpha_1,\alpha_2,\beta,\beta_1,\beta_2\}$ for which all elements have order $p$, subject to the following nontrivial commutator relations:
$$\{[\alpha_1,\alpha_2]=\beta,[\beta,\alpha_1]=\beta_1=[\alpha,\alpha_1],[\beta,\alpha_2]=\beta_2\}.$$
It follows from Lemma \ref{lem2.3} that the group $[G,G^{\vp}]$ is generated by
$[\alpha,\alpha_1^{\vp}],[\alpha_1,\alpha_2^{\vp}],[\alpha_1,\beta^{\vp}]$ and $[\alpha_2,\beta^{\vp}]$ modulo $M_0(G)$. Any two elements of these generators are commuting  modulo $M_0(G)$. In fact, by Proposition \ref{prop2.1} (4), we see that $[[\alpha,\alpha_1^{\vp}],[\alpha_1,\alpha_2^{\vp}]]=
[[\alpha,\alpha_1],[\alpha_1,\alpha_2]^{\vp}]=[\beta_1,\beta^{\vp}]\in M_0(G)$.
Similarly, for $i=1,2$, we have $[[\alpha,\alpha_1^{\vp}],[\alpha_i,\beta^{\vp}]]=[\beta_1,(\beta_i^{-1})^\vp]\in M_0(G)$, $[[\alpha_1,\alpha_2^{\vp}],[\alpha_i,\beta^{\vp}]]=[\beta,(\beta_i^{-1})^\vp]\in M_0(G)$ and
$[[\alpha_1,\beta^{\vp}],[\alpha_2,\beta^{\vp}]]=[\beta_1,(\beta_2^{-1})^\vp]\in M_0(G)$.
Thus, each element $w\in [G,G^\vp]$ can be expressed as
$$w=[\alpha,\alpha_1^{\vp}]^m[\alpha_1,\alpha_2^{\vp}]^n[\alpha_1,\beta^{\vp}]^s
[\alpha_2,\beta^{\vp}]^t\cdot w_0$$
for some $w_0\in M_0(G)$. Note that $[G,G]$ is an abelian group, it follows from Proposition \ref{prop2.2} (2) that $[G,G^\vp]$ has nilpotency class at most 2. As in the case of $\Phi_{13}$, we observe  that each $w\in M(G)$ can be written as
$w=([\alpha,\alpha_1^{\vp}][\alpha_1,\beta^{\vp}])^s\cdot w_0$
for some $w_0\in M_0(G)$. Now it suffices to show that
$[\alpha,\alpha_1^{\vp}][\alpha_1,\beta^{\vp}]\in M_0(G)$. Since
$[\alpha\alpha_1,\alpha_1\beta]=1$, it follows that
$[\alpha\alpha_1,(\alpha_1\beta)^\vp]\in M_0(G)$. Moreover,
\begin{eqnarray*}
 [\alpha\alpha_1,(\alpha_1\beta)^\vp] &=&[\alpha\alpha_1,\beta^\vp] [\alpha\alpha_1,\alpha_1^\vp] [\alpha\alpha_1,\alpha_1^\vp,\beta^\vp]\\
 &=&[\alpha,\beta^\vp][\alpha_1,\beta^\vp]
 [\alpha,\alpha_1^\vp][\beta_1,\alpha_1^\vp][\alpha,\alpha_1^\vp][\alpha\alpha_1,\alpha_1,\beta^\vp]\\
 &=& [\alpha,\beta^\vp][\alpha_1,\beta^\vp]
 [\alpha,\alpha_1^\vp][\beta_1,\alpha_1^\vp][\alpha,\alpha_1^\vp][\beta_1,\beta^\vp]
\end{eqnarray*}
which, together with the fact that $[\alpha,\alpha_1^{\vp}]$ commutes with $[\alpha_1,\beta^{\vp}]$ modulo $M_0(G)$, implies that $[\alpha,\alpha_1^{\vp}][\alpha_1,\beta^{\vp}]\in M_0(G)$. Hence, $B_{0}(\Phi_{19})=0$.
\end{proof}

\begin{prop}
$B_{0}(\Phi_{22})=0$.
\end{prop}

\begin{proof}
We take $G=\Phi_{22}(1^6)$. The polycyclic presentation of $G$ in James \cite{Jam1980}
gives rise to another polycyclic presentation of $G$ in which all generators $\alpha,\alpha_1,\alpha_2,\alpha_3,\beta_1,\beta_2$ are of order $p$ and
have the following nontrivial commutator relations:
$$[\alpha_1,\alpha]=\alpha_2,[\alpha_2,\alpha]=\alpha_3=[\beta_1,\beta_2].$$
By Lemma \ref{lem2.3}, the group $[G,G^\vp]$ is generated by $[\alpha_1,\alpha^\vp],[\alpha_2,\alpha^\vp]$ and $[\beta_1,\beta_2^\vp]$ modulo $M_0(G)$.
Note that $[[\alpha_1,\alpha^\vp],[\alpha_2,\alpha^\vp]]=
[[\alpha_1,\alpha^\vp],[\beta_1,\beta_2^\vp]]=[\alpha_2,\alpha_3^\vp]\in M_0(G)$ and
$[[\alpha_2,\alpha^\vp],[\beta_1,\beta_2^\vp]]=[\alpha_3,\alpha_3^\vp]\in M_0(G)$.
This means that every element in $[G,G^\vp]$ can be expressed as
$$[\alpha_1,\alpha^\vp]^r[\alpha_2,\alpha^\vp]^s[\beta_1,\beta_2^\vp]^t\cdot w_0$$
for some $w_0\in M_0(G)$. As $[G,G]$ is an abelian group, it follows from Proposition \ref{prop2.2} (2) that $[G,G^\vp]$ has nilpotency class at most 2. As in the case of $\Phi_{13}$, we observe  that each $w\in M(G)$ can be written as
$$w=([\alpha_2,\alpha^\vp][\beta_1,\beta_2^\vp]^{-1})^t\cdot w_0$$
for some $w_0\in M_0(G)$. To complete the proof, we need to prove that $[\alpha_2,\alpha^\vp][\beta_1,\beta_2^\vp]^{-1}\in M_0(G).$
As $[\alpha_2\beta_2,\alpha\beta_1]=1$, it follows that $[\alpha_2\beta_2,(\alpha\beta_1)^\vp]\in M_0(G)$. Extending $[\alpha_2\beta_2,(\alpha\beta_1)^\vp]$, we see that
$[\alpha_2\beta_2,(\alpha\beta_1)^\vp]=[\alpha_2,\beta_1^\vp]
[\beta_2,\beta_1^\vp][\alpha_2,\alpha^\vp][\alpha_3,\beta_2^\vp][\beta_2,\alpha^\vp]
[\alpha_3,\beta_1^\vp]$.
Note that $$[\alpha_2,\beta_1^\vp], [\alpha_3,\beta_2^\vp], [\beta_2,\alpha^\vp],
[\alpha_3,\beta_1^\vp]\in M_0(G).$$
Hence, $[\alpha_2,\alpha^\vp][\beta_1,\beta_2^\vp]^{-1}=[\alpha_2,\alpha^\vp][\beta_2,\beta_1^\vp]\in M_0(G)$.
\end{proof}

\begin{prop}
$B_{0}(\Phi_{23})=0$.
\end{prop}

\begin{proof}
We take $G=\Phi_{23}(1^6)$. The polycyclic presentation of $G$ in \cite{Jam1980}
consists of 6 generators $\alpha,\alpha_1,\alpha_2,\alpha_3,\alpha_4,\gamma$ of order $p$ and
4 nontrivial commutator relations:
$$[\alpha_1,\alpha]=\alpha_2,[\alpha_2,\alpha]=\alpha_3,[\alpha_3,\alpha]=\alpha_4,
[\alpha_1,\alpha_2]=\gamma.$$
It follows from Lemma \ref{lem2.3} that the group $[G,G^\vp]$ is generated by $[\alpha_1,\alpha^\vp],[\alpha_2,\alpha^\vp],[\alpha_3,\alpha^\vp]$ and $[\alpha_1,\alpha_2^\vp]$ modulo $M_0(G)$. Any two of these generators are commuting  modulo $M_0(G)$. Thus each element $w$ of $[G,G^\vp]$ can be expressed as
$$w=\left(\prod_{i=1}^3 [\alpha_i,\alpha^\vp]^{m_i}\right)\cdot[\alpha_1,\alpha_2^\vp]^n\cdot w_0$$
where $w_0\in M_0(G)$. Moreover, as $w\in M(G)$, we see that
$\alpha_2^{m_1}\alpha_3^{m_2}\alpha_4^{m_3}\gamma^n=1$. Hence, $m_1,m_2,m_3$ and $n$ are all divisible by $p.$

Recall that $G$ has nilpotency class 4, it follows from Proposition \ref{prop2.2} (1) that
$\tau(G)$ is nilpotent of class at most 5. Applying Proposition \ref{prop2.2} (5) on the case
$(x,y,n)=(\alpha_1,\alpha_2^\vp,p)$ in $H=\tau(G)$, we have
$$[\alpha_1^{p},\alpha_2^\vp]=[\alpha_1,\alpha_2^\vp]^{p}[\alpha_1,\alpha_2^\vp,\alpha_1]^{p\choose 2}[\alpha_1,\alpha_2^\vp,_{2}\alpha_1]^{p\choose 3}[\alpha_1,\alpha_2^\vp,_{3}\alpha_1]^{p\choose 4}[\alpha_1,\alpha_2^\vp,\alpha_1,[\alpha_1,\alpha_2^\vp]]^{\delta(p)}$$
where $\delta(p)=\frac{p(p-1)(2p-1)}{6}$. Note that
$[\alpha_1,\alpha_2^\vp,\alpha_1]^{p\choose 2}=[\alpha_1,\alpha_2,\alpha_1^\vp]^{p\choose 2}
=[\gamma,\alpha_1^\vp]^{p\choose 2}$. As $[\gamma,\alpha_1]=1$, it follows from Proposition \ref{prop2.1} (5) that $[\gamma,\alpha_1^\vp]^{p\choose 2}=[\gamma^{p\choose 2},\alpha_1^\vp]=[1,\alpha_1^\vp]=1$. We further observe that
$[\alpha_1,\alpha_2^\vp,_{2}\alpha_1]^{p\choose 3}=[\alpha_1,\alpha_2^\vp,_{3}\alpha_1]^{p\choose 4}=[\alpha_1,\alpha_2^\vp,\alpha_1,[\alpha_1,\alpha_2^\vp]]^{\delta(p)}=1.$ Thus
$$[\alpha_1,\alpha_2^\vp]^{p}=[\alpha_1^{p},\alpha_2^\vp]=[1,\alpha_2^\vp]=1.$$
Similarly, one shows that $[\alpha_1,\alpha^\vp]^{p}=[\alpha_2,\alpha^\vp]^{p}=[\alpha_3,\alpha^\vp]^{p}=1$.
Therefore, $w=w_0\in M_0(G)$ and $B_{0}(\Phi_{23})=0$.
\end{proof}

\begin{prop}
$B_{0}(\Phi_{24})=0$.
\end{prop}

\begin{proof}
Let $G=\Phi_{24}(1^6)$. The polycyclic presentation of $G$ in James \cite{Jam1980}
consists of 6 generators $\alpha,\alpha_1,\alpha_2,\alpha_3,\alpha_4,\beta$ of order $p$ and
4 nontrivial commutator relations:
$$[\alpha_1,\alpha]=\alpha_2,[\alpha_2,\alpha]=\alpha_3,[\alpha_3,\alpha]=\alpha_4=
[\alpha_1,\beta].$$
By Lemma \ref{lem2.3}, the group $[G,G^\vp]$ is generated by $[\alpha_1,\alpha^\vp],[\alpha_2,\alpha^\vp],[\alpha_3,\alpha^\vp]$ and $[\alpha_1,\beta^\vp]$ modulo $M_0(G)$. Any two of these generators are commuting  modulo $M_0(G)$. Thus each element $w$ of $[G,G^\vp]$ can be expressed as
$$w=\left(\prod_{i=1}^3 [\alpha_i,\alpha^\vp]^{m_i}\right)\cdot[\alpha_1,\beta^\vp]^n\cdot w_0$$
where $w_0\in M_0(G)$. Moreover, as $w\in M(G)$, we see that
$\alpha_2^{m_1}\alpha_3^{m_2}\alpha_4^{m_3+n}=1$. Hence, $m_1,m_2$ and $m_3+n$ are all divisible by $p.$ As in the case of $\Phi_{23}$, one can show that
$[\alpha_i,\alpha^\vp]^{p}=1=[\alpha_1,\beta^\vp]^p$ for $i=1,2,3$. Thus
$$w=([\alpha_3,\alpha^\vp][\alpha_1,\beta^\vp]^{-1})^{m_3}\cdot w_0$$
for some $w_0\in M_0(G)$. Note that $[\alpha_3\beta,\alpha\alpha_1]=1$, so $[\alpha_3\beta,(\alpha\alpha_1)^\vp]\in M_0(G)$. On the other hand,
$$[\alpha_3\beta,(\alpha\alpha_1)^\vp]=[\alpha_3,\alpha_1^\vp][\alpha_3,\alpha^\vp]
[\alpha_4,\alpha_1^\vp][\alpha_4,\beta^\vp][\beta,\alpha_1^\vp][\beta,\alpha^\vp].$$
Hence, $[\alpha_3,\alpha^\vp][\alpha_1,\beta^\vp]^{-1}\in M_0(G)$ and $w\in M_0(G)$, as desired.
\end{proof}

\begin{prop}
$B_{0}(\Phi_{25})=B_{0}(\Phi_{26})=0$.
\end{prop}

\begin{proof}
Here we only give the proof for the case $\Phi_{25}$, as the proof for the case $\Phi_{26}$ is almost same. We take $G=\Phi_{25}(222)$. The polycyclic presentation of $G$ in James \cite{Jam1980}
consists of 5 generators $\alpha,\alpha_1,\alpha_2,\alpha_3,\alpha_4$ and
the following nontrivial relations:
$$[\alpha_1,\alpha]=\alpha_2,[\alpha_2,\alpha]=\alpha_3,[\alpha_3,\alpha]=\alpha_4=\alpha_2^{(p)},
\alpha_1^{(p)}=\alpha_3, \alpha^{p^2}=\alpha_3^p=\alpha_4^p=1.$$
By Lemma \ref{lem2.3}, the group $[G,G^\vp]$ is generated by $[\alpha_1,\alpha^\vp],[\alpha_2,\alpha^\vp]$ and $[\alpha_3,\alpha^\vp]$ modulo $M_0(G)$.
For $1\leq i<j\leq 3$, we have
$$[[\alpha_i,\alpha^\vp],[\alpha_j,\alpha^\vp]]=[[\alpha_i,\alpha],[\alpha_j,\alpha]^\vp]=
[\alpha_{i+1},\alpha_{j+1}^\vp]\in M_0(G)$$
which means that $[\alpha_i,\alpha^\vp]$ and $[\alpha_j,\alpha^\vp]$ commutates modulo $M_0(G)$. Thus each element $w$ of $[G,G^\vp]$ can be expressed as
$$w=\prod_{i=1}^3 [\alpha_i,\alpha^\vp]^{m_i}\cdot w_0$$
where $w_0\in M_0(G)$. Since $w\in M(G)$, it follows that
$\alpha_2^{m_1}\alpha_3^{m_2}\alpha_4^{m_3}=1$. Hence,  $m_2$ and $m_3$ both are divisible by $p$, and $p^2$ divides $m_1$. As in the case of $\Phi_{23}$, one can show that
$[\alpha_1,\alpha^\vp]^{p^2}=[\alpha_2,\alpha^\vp]^{p}=[\alpha_3,\alpha^\vp]^p=1$.
Hence, $w=w_0\in M_0(G)$ and $B_{0}(\Phi_{25})=0$.
\end{proof}

\begin{prop}
$B_{0}(\Phi_{27})=0$.
\end{prop}

\begin{proof} Let $G=\Phi_{27}(1^6)$.
The polycyclic presentation of $G$ in James \cite{Jam1980}
consists of 6 generators $\alpha,\alpha_1,\alpha_2,\alpha_3,\alpha_4,\beta$ of order $p$ and
5 nontrivial commutator relations:
$$[\alpha_1,\alpha]=\alpha_2,[\alpha_2,\alpha]=\alpha_3,[\alpha_3,\alpha]=\alpha_4=
[\alpha_1,\beta]=[\alpha_1,\alpha_2].$$
By Lemma \ref{lem2.3}, the group $[G,G^\vp]$ is generated by $[\alpha_1,\alpha^\vp],[\alpha_2,\alpha^\vp],[\alpha_3,\alpha^\vp],[\alpha_1,\beta^\vp]$ and $[\alpha_1,\alpha_2^\vp]$ modulo $M_0(G)$. One can check that any two of these generators are commuting  modulo $M_0(G)$. Thus each element $w$ of $[G,G^\vp]$ can be expressed as
$$w=\left(\prod_{i=1}^3 [\alpha_i,\alpha^\vp]^{m_i}\right)\cdot[\alpha_1,\beta^\vp]^s
\cdot[\alpha_1,\alpha_2^\vp]^t\cdot w_0$$
where $w_0\in M_0(G)$. As $w\in M(G)$, we see that $\alpha_2^{m_1}\alpha_3^{m_2}\alpha_4^{m_3+s+t}=1$. Hence, $m_1,m_2$ and $m_3+s+t$ are all divisible by $p.$ As in the case of $\Phi_{23}$, one can show that
$[\alpha_1,\alpha_2^\vp]^p=[\alpha_1,\beta^\vp]^p=[\alpha_i,\alpha^\vp]^p=1$ where $i=1,2,3$.
Thus,
$$w=[\alpha_3,\alpha^\vp]^{m_3}\cdot[\alpha_1,\beta^\vp]^{-m_3-t}
\cdot[\alpha_1,\alpha_2^\vp]^t\cdot w_0.$$
To complete the proof, it suffices to show that $[\alpha_3,\alpha^\vp]
[\alpha_1,\beta^\vp]^{-1}$ and $[\alpha_1,\beta^\vp]^{-1}[\alpha_1,\alpha_2^\vp]$
both are in $M_0(G)$. The fact that $[\alpha_3\beta,\alpha\alpha_1]=1$ implies that
$[\alpha_3\beta,(\alpha\alpha_1)^\vp]\in M_0(G)$. Note that
$$[\alpha_3\beta,(\alpha\alpha_1)^\vp]=[\alpha_3,\alpha_1^\vp][\alpha_3,\alpha^\vp]
[\alpha_4,\alpha_1^\vp][\alpha_4,\beta^\vp][\beta,\alpha_1^\vp][\beta,\alpha^\vp]$$
and $[\alpha_3,\alpha_1^\vp],
[\alpha_4,\alpha_1^\vp],[\alpha_4,\beta^\vp],[\beta,\alpha^\vp]$ belong to $M_0(G)$. Hence,
$[\alpha_3,\alpha^\vp][\alpha_1,\beta^\vp]^{-1}\in M_0(G)$. Similarly, as
$[\beta\alpha_1,\alpha_1\alpha_2]=1$, one can expand $[\beta\alpha_1,(\alpha_1\alpha_2)^\vp]$ to see that $[\alpha_1,\beta^\vp]^{-1}[\alpha_1,\alpha_2^\vp]\in M_0(G)$.
Therefore, $w\in M_0(G)$ and $B_{0}(\Phi_{27})=0$.
\end{proof}

\begin{prop}
$B_{0}(\Phi_{30})=0$.
\end{prop}

\begin{proof}
We take $G=\Phi_{30}(1^6)$. The polycyclic presentation of $G$ in James \cite{Jam1980}
consists of 6 generators $\alpha,\alpha_1,\alpha_2,\alpha_3,\alpha_4,\beta$ of order $p$ and
5 nontrivial commutator relations:
$$[\alpha_1,\alpha]=\alpha_2,[\alpha_2,\alpha]=\alpha_3=[\alpha_1,\beta],[\alpha_2,\beta]=\alpha_4
=[\alpha_3,\alpha].$$
By Lemma \ref{lem2.3}, the group $[G,G^\vp]$ is generated by $[\alpha_1,\alpha^\vp],[\alpha_2,\alpha^\vp],[\alpha_1,\beta^\vp],[\alpha_2,\beta^\vp]$ and $[\alpha_3,\alpha^\vp]$ modulo $M_0(G)$. As before, one can check that any two of these generators are commuting  modulo $M_0(G)$. Thus each element $w$ of $[G,G^\vp]$ can be expressed as
$$w=[\alpha_1,\alpha^\vp]^m[\alpha_2,\alpha^\vp]^n[\alpha_1,\beta^\vp]^r[\alpha_2,\beta^\vp]^s
[\alpha_3,\alpha^\vp]^t\cdot w_0$$
where $w_0\in M_0(G)$. As $w\in M(G)$, we see that $\alpha_2^{m}\alpha_3^{n+r}\alpha_4^{s+t}=1$. Hence, $m,n+r$ and $s+t$ are all divisible by $p.$ As in the case of $\Phi_{23}$, one can show that
$$[\alpha_1,\alpha^\vp]^p=[\alpha_2,\alpha^\vp]^p=[\alpha_1,\beta^\vp]^p=[\alpha_2,\beta^\vp]^p=
[\alpha_3,\alpha^\vp]^p=1.$$
Thus,
$$w=([\alpha_2,\alpha^\vp][\alpha_1,\beta^\vp]^{-1})^n\cdot
([\alpha_2,\beta^\vp][\alpha_3,\alpha^\vp]^{-1})^s\cdot w_0.$$
To complete the proof, it suffices to show that $[\alpha_2,\alpha^\vp][\alpha_1,\beta^\vp]^{-1}$ and $[\alpha_2,\beta^\vp][\alpha_3,\alpha^\vp]^{-1}$
both belong to $M_0(G)$. Indeed, the fact that $[\alpha_2\beta,\alpha\alpha_1]=1$ implies that
$[\alpha_2\beta,(\alpha\alpha_1)^\vp]\in M_0(G)$. On the other hand,
$$[\alpha_2\beta,(\alpha\alpha_1)^\vp]=[\alpha_2,\alpha_1^\vp][\alpha_2,\alpha^\vp]
[\alpha_3,\alpha_1^\vp][\alpha_3,\beta^\vp][\beta,\alpha_1^\vp][\beta,\alpha^\vp]$$
and $[\alpha_2,\alpha_1^\vp],
[\alpha_3,\alpha_1^\vp],[\alpha_3,\beta^\vp],[\beta,\alpha^\vp]$ all belong to $M_0(G)$. Hence,
$[\alpha_2,\alpha^\vp][\alpha_1,\beta^\vp]^{-1}\in M_0(G)$. Similarly, as
$[\alpha_2\alpha,\beta\alpha_3]=1$, one can expand $[\alpha_2\alpha,(\beta\alpha_3)^\vp]$ to see that $[\alpha_2,\beta^\vp][\alpha_3,\alpha^\vp]^{-1}\in M_0(G)$.
Therefore, $w\in M_0(G)$ and $B_{0}(\Phi_{30})=0$.
\end{proof}

\begin{prop}
$B_{0}(\Phi_{31})=B_{0}(\Phi_{33})=B_{0}(\Phi_{34})=0$.
\end{prop}

\begin{proof}
We take $G=\Phi_{31}(1^6)$. The polycyclic presentation of $G$ in James \cite{Jam1980}
consists of 6 generators $\alpha,\alpha_1,\alpha_2,\beta_1,\beta_2,\gamma$ of order $p$ and
4 nontrivial commutator relations:
$$[\alpha_1,\alpha]=\beta_1,[\alpha_2,\alpha]=\beta_2,[\alpha_1,\beta_1]=\gamma=
[\alpha_2,\beta_2].$$
By Lemma \ref{lem2.3}, the group $[G,G^\vp]$ is generated by $[\alpha_1,\alpha^\vp],[\alpha_2,\alpha^\vp],[\alpha_1,\beta_1^\vp]$ and $[\alpha_2,\beta_2^\vp]$ modulo $M_0(G)$. As before, one can check that any two of these generators are commuting  modulo $M_0(G)$. Thus each element $w$ of $[G,G^\vp]$ can be expressed as
$$w=[\alpha_1,\alpha^\vp]^m[\alpha_2,\alpha^\vp]^n[\alpha_1,\beta_1^\vp]^s[\alpha_2,\beta_2^\vp]^t
\cdot w_0$$
where $w_0\in M_0(G)$. As $w\in M(G)$, we see that $\beta_1^{m}\beta_2^{n}\gamma^{s+t}=1$. Hence, $m,n$ and $s+t$ are all divisible by $p.$ As in the case of $\Phi_{13}$, one can show that
$[\alpha_1,\alpha^\vp]^p=[\alpha_2,\alpha^\vp]^p=[\alpha_1,\beta_1^\vp]^p=[\alpha_2,\beta_2^\vp]^p
=1.$
Thus,
$$w=([\alpha_1,\beta_1^\vp][\alpha_2,\beta_2^\vp]^{-1})^s\cdot w_0.$$
It is sufficient to show that $[\alpha_1,\beta_1^\vp][\alpha_2,\beta_2^\vp]^{-1}\in M_0(G)$. Indeed, the fact that $[\alpha_1\beta_2,\beta_1\alpha_2]=1$ implies that
$[\alpha_1\beta_2,(\beta_1\alpha_2)^\vp]\in M_0(G)$. On the other hand,
$$[\alpha_1\beta_2,(\beta_1\alpha_2)^\vp]=[\alpha_1,\alpha_2^\vp][\beta_2,\alpha_2^\vp]
[\alpha_1,\beta_1^\vp][\gamma,\beta_2^\vp][\beta_2,\beta_1^\vp][\gamma,\alpha_2^\vp]$$
and $[\alpha_1,\alpha_2^\vp],
[\gamma,\beta_2^\vp],[\gamma,\alpha_2^\vp],[\beta_2,\beta_1^\vp]$ all belong to $M_0(G)$. Hence,
$[\alpha_1,\beta_1^\vp][\alpha_2,\beta_2^\vp]^{-1}\in M_0(G)$.
Therefore, $w\in M_0(G)$ and $B_{0}(\Phi_{31})=0$.

Similar argument applies on the cases of $\Phi_{33}$ and $\Phi_{34}$. Here we only sketch the proofs. First, we take $G=\Phi_{33}(1^6)$ so that it has a polycyclic presentation of $G$ in James \cite{Jam1980}
consists of 6 generators $\alpha,\alpha_1,\alpha_2,\beta_1,\beta_2,\gamma$ of order $p$ and
4 nontrivial commutator relations:
$[\alpha_1,\alpha]=\beta_1,[\alpha_2,\alpha]=\beta_2,[\alpha_1,\beta_1]=\gamma=
[\beta_2,\alpha].$  Note that each element $w\in [G,G^\vp]$
can be expressed as
$$w=[\alpha_1,\alpha^\vp]^m[\alpha_2,\alpha^\vp]^n[\alpha_1,\beta_1^\vp]^s[\beta_2,\alpha^\vp]^t
\cdot w_0$$
for some $w_0\in M_0(G)$. Moreover, as $w\in M(G)$, one can show that
$w=([\alpha_1,\beta_1^\vp][\beta_2,\alpha^\vp]^{-1})^s\cdot w_0$. Note that $[\alpha_1\alpha,(\beta_1\beta_2)^\vp]\in M_0(G)$ and
$$[\alpha_1\alpha,(\beta_1\beta_2)^\vp]=[\alpha_1,\alpha_2^\vp][\alpha,\beta_2^\vp]
[\alpha_1,\beta_1^\vp][\gamma,\alpha^\vp][\alpha,\beta_1^\vp][\gamma,\beta_2^\vp].$$
Since $[\alpha_1,\alpha_2^\vp],[\gamma,\alpha^\vp],[\alpha,\beta_1^\vp],[\gamma,\beta_2^\vp]$ all belong to $M_0(G)$, it follows that $[\alpha_1,\beta_1^\vp][\beta_2,\alpha^\vp]^{-1}\in M_0(G)$. Hence, $M(G)=M_0(G)$ and $B_{0}(\Phi_{33})=0$.

For the case of $\Phi_{34}$, we take $G=\Phi_{34}(321)$. The polycyclic presentation of $G$ in James \cite{Jam1980}
consists of 6 generators $\alpha,\alpha_1,\alpha_2,\beta_1,\beta_2,\gamma$  and
the following nontrivial relations:
$$[\alpha_1,\alpha]=\beta_1,[\alpha_2,\alpha]=\beta_2,[\alpha_1,\beta_1]=\gamma=\beta_1^p=
[\beta_2,\alpha],\alpha^p=\beta_1,\alpha_1^p=\beta_2.$$
Note that each element $w\in [G,G^\vp]$
can be expressed as
$$w=[\alpha_1,\alpha^\vp]^m[\alpha_2,\alpha^\vp]^n[\beta_2,\alpha^\vp]^s[\alpha_1,\beta_1^\vp]^t
\cdot w_0$$
for some $w_0\in M_0(G)$. Moreover, as $w\in M(G)$, one can show that $n$ and $s+t$ both are divisible by $p$, and $p^2$ divides $m$. One can show that 
$$[\alpha_1,\alpha^\vp]^{p^2}=[\alpha_2,\alpha^\vp]^p=[\beta_2,\alpha^\vp]^p=
[\alpha_1,\beta_1^\vp]^p=1.$$ Thus
$w=([\beta_2,\alpha^\vp][\alpha_1,\beta_1^\vp]^{-1})^t\cdot w_0$.
Since
$[\beta_2\beta_1,\alpha\alpha_1]=1$, one can expand $[\beta_2\beta_1,(\alpha\alpha_1)^\vp]$ to see that $[\beta_2,\alpha^\vp][\alpha_1,\beta_1^\vp]^{-1}\in M_0(G)$.
Hence, $M(G)=M_0(G)$ and $B_{0}(\Phi_{34})=0$.
\end{proof}

\begin{prop}
$B_{0}(\Phi_{35})=B_{0}(\Phi_{37})=0$.
\end{prop}

\begin{proof}
We take $G=\Phi_{35}(1^6)$ and $\Phi_{37}(1^6)$ respectively. They both are generated by 
$\alpha,\alpha_1,\dots,\alpha_5$. 

For the case of $\Phi_{35}$, the generators have the following nontrivial commutator relations:
$$[\alpha_i,\alpha]=\alpha_{i+1}$$
for $1\leq i\leq 4.$ 
By Lemma \ref{lem2.3}, the group $[G,G^\vp]$ is generated by $\{[\alpha_i,\alpha^\vp]\mid 1\leq i\leq 4\}$ modulo $M_0(G)$. As before, one can check that any two of these generators of $[G,G^\vp]$ are commuting  modulo $M_0(G)$. Thus each element $w$ of $[G,G^\vp]$ can be expressed as
$$w=\left(\prod_{i=1}^4[\alpha_i,\alpha^\vp]^{m_i}\right)\cdot w_0$$
where $w_0\in M_0(G)$. As $w\in M(G)$, we see that $\alpha_2^{m_1}\alpha_3^{m_2}\alpha_4^{m_3}\alpha_5^{m_4}=1$. Hence, $p$ divides $m_i$ for $1\leq i\leq 4.$  As in the case of $\Phi_{23}$, one can show that
$[\alpha_i,\alpha^\vp]^p=1$ for $1\leq i\leq 4.$
Hence, $w=w_0\in M_0(G)$ and $B_{0}(\Phi_{35})=0$.

The case of $\Phi_{37}$ is much more complicated. Firstly, the generators for $G=\Phi_{37}(1^6)$ have the following nontrivial commutator relations:
$$[\alpha_i,\alpha]=\alpha_{i+1}, [\alpha_2,\alpha_3]=[\alpha_3,\alpha_1]=[\alpha_4,\alpha_1]=\alpha_5$$
for $i=1,2,3.$ The group $[G,G^\vp]$ can be generated by 
$$[\alpha_1,\alpha^\vp],[\alpha_2,\alpha^\vp],[\alpha_3,\alpha^\vp],
[\alpha_2,\alpha_3^\vp],[\alpha_3,\alpha_1^\vp],[\alpha_4,\alpha_1^\vp]$$
modulo $M_0(G)$. 
Except for $[\alpha_1,\alpha^\vp]$ and $[\alpha_2,\alpha^\vp]$, 
any two of these generators commutates each other modulo $M_0(G)$.  
Note that $[[\alpha_1,\alpha^\vp],[\alpha_2,\alpha^\vp]]=[\alpha_2,\alpha_3^\vp]$. Thus 
each element $w\in [G,G^\vp]$
can be expressed as
$$w=\left(\prod_{i=1}^3[\alpha_i,\alpha^\vp]^{m_i}\right)\cdot[\alpha_2,\alpha_3^\vp]^r\cdot
[\alpha_3,\alpha_1^\vp]^s\cdot[\alpha_4,\alpha_1^\vp]^t
\cdot w_0$$
for some $w_0\in M_0(G)$. Moreover, as $w\in M(G)$, we see that $p$ divides $m_1,m_2,m_3$ and 
$r+s+t$ respectively. One can show that
$[\alpha_1,\alpha^\vp]^p=[\alpha_2,\alpha^\vp]^p=[\alpha_3,\alpha^\vp]^p=
[\alpha_2,\alpha_3^\vp]^p=[\alpha_3,\alpha_1^\vp]^p=[\alpha_4,\alpha_1^\vp]^p=1.$ Thus
$$w=([\alpha_2,\alpha_3^\vp][\alpha_3,\alpha_1^\vp]^{-1})^r
([\alpha_3,\alpha_1^\vp]^{-1}[\alpha_4,\alpha_1^\vp])^t\cdot w_0.$$
We observe that $[\alpha_2\alpha_1,\alpha_3]=[\alpha_2,\alpha_3][\alpha_2,\alpha_3,\alpha_1][\alpha_1,\alpha_3]=1$.
Thus $[\alpha_2\alpha_1,\alpha_3^\vp]\in M_0(G)$.
Note that $[\alpha_2\alpha_1,\alpha_3^\vp]=
[\alpha_2,\alpha_3^\vp][\alpha_2,\alpha_3^\vp,\alpha_1][\alpha_1,\alpha_3^\vp].$
As $[\alpha_2,\alpha_3^\vp,\alpha_1]=[\alpha_5,\alpha_1^\vp]\in M_0(G)$, 
it follows that $[\alpha_2,\alpha_3^\vp][\alpha_3,\alpha_1^\vp]^{-1}=
[\alpha_2,\alpha_3^\vp][\alpha_1,\alpha_3^\vp]\in M_0(G)$.
Moreover, as $[\alpha_3\alpha_1,\alpha_1\alpha_4]=1$, it follows that
$[\alpha_3\alpha_1,(\alpha_1\alpha_4)^\vp]\in M_0(G)$. Expanding 
$[\alpha_3\alpha_1,(\alpha_1\alpha_4)^\vp]$, we see that 
$[\alpha_3,\alpha_1^\vp]^{-1}[\alpha_4,\alpha_1^\vp]\in   M_0(G)$.
Therefore, $w\in M_0(G)$ and $B_{0}(\Phi_{37})=0$.
\end{proof}

\begin{prop}
$B_{0}(\Phi_{40})=B_{0}(\Phi_{41})=B_{0}(\Phi_{42})=B_{0}(\Phi_{43})=0$.
\end{prop}

\begin{proof}
We will take $G=\Phi_{40}(1^6),\Phi_{41}(1^6),\Phi_{42}(222)_{a_0}$ and $\Phi_{43}(222)_{a_r}$ respectively. They both are generated by
$\alpha_1,\alpha_2,\beta,\beta_1,\beta_2,\gamma$.

(1) We take $G=\Phi_{40}(1^6)$. The $G$ has the following nontrivial commutator relations:
$$[\alpha_1,\alpha_2]=\beta,[\beta,\alpha_1]=\beta_1,[\beta,\alpha_2]=\beta_2,
[\beta_1,\alpha_2]=[\beta_2,\alpha_1]=\gamma.$$
By Lemma \ref{lem2.3}, the group $[G,G^\vp]$ is generated by $[\alpha_1,\alpha_2^\vp],
[\beta,\alpha_1^\vp], [\beta,\alpha_2^\vp], [\beta_1,\alpha_2^\vp]$ and $[\beta_2,\alpha_1^\vp]$ modulo $M_0(G)$. One can check that any two of these generators of $[G,G^\vp]$ are commuting  modulo $M_0(G)$. Thus each element $w$ of $[G,G^\vp]$ can be expressed as
$$w=[\alpha_1,\alpha_2^\vp]^m
[\beta,\alpha_1^\vp]^n[\beta,\alpha_2^\vp]^r [\beta_1,\alpha_2^\vp]^s [\beta_2,\alpha_1^\vp]^t\cdot w_0$$
where $w_0\in M_0(G)$. As $w\in M(G)$, we see that $\beta^{m}\beta_1^{n}\beta_2^{r}\gamma^{s+t}=1$. Hence, $p$ divides $m,n,r$ and $s+t$ respectively.  As in the case of $\Phi_{23}$, one can show that
$$[\alpha_1,\alpha_2^\vp]^p=
[\beta,\alpha_1^\vp]^p=[\beta,\alpha_2^\vp]^p=[\beta_1,\alpha_2^\vp]^p= [\beta_2,\alpha_1^\vp]^p=1$$
which implies that 
$$w=([\beta_1,\alpha_2^\vp][\beta_2,\alpha_1^\vp]^{-1})^{s}\cdot w_0.$$
Note that $[\alpha_1\alpha_2\beta_1,\alpha_1\alpha_2\beta_2]=1$. Thus
$[\alpha_1\alpha_2\beta_1,(\alpha_1\alpha_2\beta_2)^\vp]\in M_0(G)$. On the other hand,
\begin{eqnarray*}
[\alpha_1\alpha_2\beta_1,(\alpha_1\alpha_2\beta_2)^\vp]&=&
[\alpha_1,\beta_2^\vp][\alpha_1,\beta_2^\vp,\alpha_2][\alpha_2,\beta_2^\vp]
[\gamma^{-1},\beta_1^\vp][\beta_1,\beta_2^\vp]\cdot\\
&&[\beta_1,\alpha_2^\vp][\beta_1,\alpha_2^\vp,\alpha_1^\vp]
[\beta_1,\alpha_1^\vp][\gamma,\beta_2^\vp].
\end{eqnarray*}
In the two-hand sides of the previous equality, all commutator factors except for $[\beta_1,\alpha_2^\vp]$ and $[\alpha_1,\beta_2^\vp]$, are in $M_0(G)$.
Hence, $[\beta_1,\alpha_2^\vp][\beta_2,\alpha_1^\vp]^{-1}\in M_0(G)$ and $B_{0}(\Phi_{40})=0.$

(2) We take $G=\Phi_{41}(1^6)$. The $G$ has the following nontrivial commutator relations:
$$[\alpha_1,\alpha_2]=\beta,[\beta,\alpha_1]=\beta_1,[\beta,\alpha_2]=\beta_2,
[\alpha_1,\beta_1]^{-\nu}=[\alpha_2,\beta_2]=\gamma^{-\nu}$$
where all generators have order $p$.
As before, one can show that each element $w\in M(G)$ can be expressed as 
$$w=([\alpha_1,\beta_1^\vp]^\nu[\alpha_2,\beta_2^\vp])^t\cdot w_0$$
for some $w_0\in M_0(G)$ and $t\in\N$. Note that $[\alpha_1\alpha_2,\beta_1^{\nu}\beta_2]=1$. Thus 
$[\alpha_1\alpha_2,(\beta_1^{\nu}\beta_2)^\vp]\in M_0(G)$. We expand $[\alpha_1\alpha_2,(\beta_1^{\nu}\beta_2)^\vp]$ as follows:
$$[\alpha_1\alpha_2,(\beta_1^{\nu}\beta_2)^\vp]=[\alpha_1,\beta_2^\vp]
[\alpha_1,\beta_1^\vp]^\nu[\gamma^\nu,\beta_2^\vp][\gamma^\nu,\alpha_2^\vp][\alpha_2,\beta_2^\vp]
[\alpha_2,\beta_1^\vp]^\nu.$$
As $[\alpha_1,\beta_2^\vp], [\gamma^\nu,\beta_2^\vp],[\gamma^\nu,\alpha_2^\vp]$ and $
[\alpha_2,\beta_1^\vp]$ all belong to $M_0(G)$, we see that 
$[\alpha_1,\beta_1^\vp]^\nu[\alpha_2,\beta_2^\vp]\in M_0(G)$. 
Hence, $w\in M_0(G)$ and $B_{0}(\Phi_{41})=0.$

(3) The group $G=\Phi_{42}(222)_{a_{0}}$  has the following nontrivial commutator relations:
$$[\alpha_1,\alpha_2]=\beta,[\beta,\alpha_1]=\beta_1,[\beta,\alpha_2]=\beta_2,
[\alpha_1,\beta_2]=[\alpha_2,\beta_1]=\beta^{p}=\gamma$$
together with $\alpha_{1}^{p}=\beta_{1}^{-1}\gamma^{-1/2},\alpha_{2}^{p}=\beta_{2}\gamma^{1/2}$ and $\beta_{1}^{p}=\beta_{2}^{p}=\gamma^{p}=1$.
One can show that each element $w\in M(G)$ can be expressed as 
$$w=([\alpha_1,\beta_2^\vp][\alpha_2,\beta_1^\vp]^{-1})^s\cdot w_0$$
for some $w_0\in M_0(G)$ and $s\in\N$. Note that $[\alpha_1\alpha_2\beta_{1},\alpha_{1}\alpha_{2}\beta_2]=1$. Thus 
$[\alpha_1\alpha_2\beta_{1},(\alpha_{1}\alpha_{2}\beta_2)^\vp]\in M_0(G)$. Expanding 
 $[\alpha_1\alpha_2\beta_{1},(\alpha_{1}\alpha_{2}\beta_2)^\vp]$, we have
$$[\alpha_1\alpha_2\beta_{1},(\alpha_{1}\alpha_{2}\beta_2)^\vp]=[\alpha_1,\beta_2^\vp]
[\alpha_{1},\beta_{2}^{\vp},\alpha_{2}][\alpha_{2},\beta_{2}^{\vp}][\gamma,\beta_{1}^{\vp}][\beta_{1},\beta_{2}^{\vp}][\beta_{1},\alpha_{2}^{\vp}][\beta_{1},\alpha_{1}^{\vp}][\beta_{1},\alpha_{1}^{\vp},\alpha_{2}^{\vp}].$$
In the two-hand sides of this equality, all commutator factors except for $[\beta_1,\alpha_2^\vp]$ and $[\alpha_1,\beta_2^\vp]$, are in $M_0(G)$.
Hence, $[\alpha_1,\beta_2^\vp][\alpha_2,\beta_1^\vp]^{-1}\in M_0(G)$ and $B_{0}(\Phi_{42})=0.$

(4) The group $G=\Phi_{43}(222)_{a_{r}}$  has the following nontrivial commutator relations:
$$[\alpha_1,\alpha_2]=\beta,[\beta,\alpha_1]=\beta_1,[\beta,\alpha_2]=\beta_2,
[\alpha_1,\beta_2]=[\alpha_2,\beta_1]=\beta^{p}=\gamma$$
together with $\alpha_{1}^{p}=\beta_{2}\gamma^{k},\alpha_{2}^{p}=\beta_{1}^{\nu}\gamma^{\ell},\beta^{p}=\gamma^{n}$, and $\beta_{1}^{p}=\beta_{2}^{p}=\gamma^{p}=1$, where
$n=\nu+ {p \choose 3}$, and $k,\ell$ are the smallest positive integers satisfying $$(k-\nu)^{2}-\nu(\ell+\nu)^{2}\equiv r \mod p,$$ for 
$r\in \{0,1,\dots,p-1\}$. As before, one can show that each element $w\in M(G)$ can be expressed as 
$$w=([\alpha_1,\beta_1^\vp]^\nu[\alpha_2,\beta_2^\vp])^t\cdot w_0$$
for some $w_0\in M_0(G)$ and $t\in\N$. Note that $[\alpha_1\alpha_2,\beta_1^{\nu}\beta_2]=1$. Thus 
$[\alpha_1\alpha_2,(\beta_1^{\nu}\beta_2)^\vp]\in M_0(G)$. We expand $[\alpha_1\alpha_2,(\beta_1^{\nu}\beta_2)^\vp]$ as follows:
$$[\alpha_1\alpha_2,(\beta_1^{\nu}\beta_2)^\vp]=[\alpha_1,\beta_2^\vp]
[\alpha_1,\beta_1^\vp]^\nu[\gamma^\nu,\beta_2^\vp][\gamma^\nu,\alpha_2^\vp][\alpha_2,\beta_2^\vp]
[\alpha_2,\beta_1^\vp]^\nu.$$
As $[\alpha_1,\beta_2^\vp], [\gamma^\nu,\beta_2^\vp],[\gamma^\nu,\alpha_2^\vp]$ and $
[\alpha_2,\beta_1^\vp]$ all belong to $M_0(G)$, we see that 
$[\alpha_1,\beta_1^\vp]^\nu[\alpha_2,\beta_2^\vp]\in M_0(G)$. 
Hence, $w\in M_0(G)$ and $B_{0}(\Phi_{43})=0.$ We are done.
\end{proof}

\section{\scshape Non-vanishing  Bogomolov Multipliers}
\setcounter{equation}{0}
\renewcommand{\theequation}
{4.\arabic{equation}}
\setcounter{theorem}{0}
\renewcommand{\thetheorem}
{4.\arabic{theorem}}

In this section, we apply the following lemma, which appeared in Hoshi-Kang-Kunyavskii \cite[Lemma 2.1]{HKK2013}, to prove that
$B_{0}(\Phi_{k})$ is not zero for $k\in \{18,20,21,36,38,39\}$. Throughout we denote the cyclic group of order $n$ by $C_{n}$ and take the convention that ${\ell\choose s}=0$
for $1\leqslant \ell<s$.

\begin{lem}\label{lem4.1}
Let $N$ be a normal subgroup of a finite group $G$. Assume that the transgression map
$\tr:H^{1}(N,\Q/\Z)^{G}\ra H^{2}(G/N,\Q/\Z)$ is not surjective and the group $EN/N$
is a cyclic subgroup of $G/N$ for any bicyclic subgroup $E$ of $G$. Then $B_{0}(G)\neq 0$.
\end{lem}

\begin{thm}
$B_{0}(\Phi_{k})\neq 0$ for $k\in \Delta=\{18,20,21,36,38,39\}$. 
\end{thm}

\begin{proof}
%The statement that $B_{0}(\Phi_{10})\neq 0$ has been proved in Proposition \ref{prop3.1}. 

Here we take $G=\Phi_{k}(1^{6})$ and  only give the proof for the case
$k=18$. Similar arguments can be applied to the remaining cases if we take $N=\langle \beta,\beta_{1},\beta_{2}\rangle$ for the cases $k\in\{20,21\}$ and take $N=\langle \alpha_{3},\alpha_{4},\alpha_{5}\rangle$ for the cases $k\in\{36,38,39\}$.

The group  $G=\Phi_{18}(1^{6})$ has a polycyclic presentation which consists of 6 generators of order $p$:  $\alpha,\alpha_{1},\alpha_{2},\alpha_{3},\beta,\gamma$ and nontrivial commutator relations:
$$[\alpha_{1},\alpha]=\alpha_{2}, [\alpha_{2},\alpha]=\alpha_{3}=[\alpha_{1},\beta],[\alpha,\beta]=\gamma.$$
Let $N$ be the subgroup of $G$ generated by $\alpha_{3},\beta,\gamma$. Then $N$ is normal and isomorphic to $C_{p}\times C_{p}\times C_{p}$. Hence, $H^{1}(N,\Q/\Z)=\Hom(N,\Q/\Z)\cong C_{p}\times C_{p}\times C_{p}$. We define
$\varphi_{1},\varphi_{2},\varphi_{3}\in H^{1}(N,\Q/\Z)$ by
\begin{eqnarray*}
\varphi_{1}(\alpha_{3})=1/p, &\varphi_{1}(\beta)=0,& \varphi_{1}(\gamma)=0; \\
\varphi_{2}(\alpha_{3})=0, &\varphi_{2}(\beta)=1/p,& \varphi_{2}(\gamma)=0; \\
\varphi_{3}(\alpha_{3})=0, &\varphi_{3}(\beta)=0,& \varphi_{3}(\gamma)=1/p.
\end{eqnarray*}
Then $H^{1}(N,\Q/\Z)=\langle \varphi_{1},\varphi_{2},\varphi_{3}\rangle$. The actions of $G$ on $H^{1}(N,\Q/\Z)$ are given by
\begin{eqnarray*}
(\alpha\cdot\varphi_{1})(\alpha_{3})&=&\varphi_{1}(\alpha^{-1}\alpha_{3}\alpha)=\varphi_{1}(\alpha_{3})=1/p\\
(\alpha\cdot\varphi_{1})(\beta)&=&\varphi_{1}(\alpha^{-1}\beta\alpha)=\varphi_{1}(\beta\gamma^{-1})=
\varphi_{1}(\beta)+\varphi_{1}(\gamma^{-1})=0\\
(\alpha\cdot\varphi_{1})(\gamma)&=& \varphi_{1}(\alpha^{-1}\gamma\alpha)=\varphi_{1}(\gamma)=0.
\end{eqnarray*}
This means that $\alpha$ fixes $\varphi_{1}$. Similarly, one can show that $\alpha$ also fixes $\varphi_{2}$ and $\alpha\cdot\varphi_{3}=-\varphi_{2}+\varphi_{3}$. Moreover, $\alpha_{1}(\varphi_{1})=\varphi_{1}-\varphi_{2},
\alpha_{1}(\varphi_{2})=\varphi_{2},\alpha_{1}(\varphi_{3})=\varphi_{3}$, and $\alpha_{2}(\varphi_{i})=\varphi_{i}$ for $i=1,2,3.$
For any $\varphi\in H^{1}(N,\Q/\Z)$, we may write $\varphi=a_{1}\varphi_{1}+a_{2}\varphi_{2}+a_{3}\varphi$ for 
$a_{1},a_{2},a_{3}\in\{0,1,\dots,p-1\}$. We observe that $\varphi\in H^{1}(N,\Q/\Z)^{G}$ if and only if $a_{1}=a_{3}=0.$
Clearly, $\varphi_{2}\in H^{1}(N,\Q/\Z)^{G}$ and $H^{1}(N,\Q/\Z)^{G}=\langle \varphi_{2}\rangle\cong C_{p}$.
Since $G/N$ is a group of order $p^{3}$ and of exponent $p$, it follows from \cite[Proposition 6.3]{Lew1968} that
$H^{2}(G/N,\Q/\Z)\cong C_{p}\times C_{p}$. Hence, the transgression map
$\tr:H^{1}(N,\Q/\Z)^{G}\ra H^{2}(G/N,\Q/\Z)$ is not surjective.

To apply Lemma \ref{lem4.1} to show that $B_{0}(G)$ is not zero, it is sufficient to show that 
the group $EN/N$ is a cyclic subgroup of $G/N$ for any bicyclic subgroup $E$ of $G$. Recall that a group $E$
is bicyclic if it is either isomorphic to a cyclic group or the direct product of two cyclic groups. By the commutator relations of $G$, we see that
\begin{eqnarray}
\alpha_{2}^{i}\cdot\alpha^{j}& = & \alpha^{j}\cdot\alpha_{2}^{i}\cdot\alpha_{3}^{ij} \label{eq4.1}\\
\alpha_{1}^{i}\cdot\beta^{j}& = & \beta^{j}\cdot\alpha_{1}^{i}\cdot\alpha_{3}^{ij} \label{eq4.2}\\
\alpha_{1}^{i}\cdot\alpha^{j}& = & \alpha^{j}\cdot\alpha_{1}^{i}\cdot\alpha_{2}^{ij}\cdot\alpha_{3}^{i\cdot{j\choose 2}} \label{eq4.3}\\
\alpha^{i}\cdot \beta^{j}& = & \beta^{j}\cdot\alpha^{i}\cdot \gamma^{ij} \label{eq4.4}
\end{eqnarray}
where $1\leqslant i,j\leqslant p-1$. Let $E$ be a bicyclic subgroup of $G$ generated by $e_{1},e_{2}$. Then $EN/N$ is a proper subgroup of $G/N$. Thus the order of $EN/N$ is either $p$ or $p^{2}$.

Now we \textbf{claim} that $|EN/N|=p$. We assume by the way of contradiction that $|EN/N|=p^{2}$. In the quotient group $G/N$,
we may write $e_{1}N=\alpha^{a_{1}}\alpha_{1}^{a_{2}}\alpha_{2}^{a_{3}}N$ and $e_{2}N=\alpha^{b_{1}}\alpha_{1}^{b_{2}}\alpha_{2}^{b_{3}}N$ for some $a_{i},b_{i}\in \N$. As in the proof of \cite[Lemma 2.1, Step 2]{HKK2013}, one shows that $(e_{1}N,e_{2}N)$ only could be one of the following three possibilities:
$$\{(\alpha_{1}N,\alpha_{2}N),(\alpha\alpha_{2}^{a_{3}}N,\alpha_{1}\alpha_{2}^{b_{3}}N),(\alpha\alpha_{1}^{a_{2}}N,\alpha_{2}N)\}$$ by changing suitable generators $e_{1},e_{2}$ if necessary. For the first case, we write $e_{1}=\alpha_{1}\alpha_{3}^{a_{4}}\beta^{a_{5}}\gamma^{a_{6}}$ and $e_{2}=\alpha_{2}\alpha_{3}^{b_{4}}\beta^{b_{5}}\gamma^{b_{6}}$. As $e_{1}$ and $e_{2}$ are commutating,  it follows that $\alpha_{1}\beta^{a_{5}}\alpha_{2}\beta^{b_{5}}=\alpha_{2}\beta^{b_{5}}\alpha_{1}\beta^{a_{5}}$, which together with (\ref{eq4.2}) implies that $[\alpha_{1},\alpha_{2}]\neq 1$. This is a contradiction. 
For the second case, using (\ref{eq4.1}) and (\ref{eq4.3}) we see that $e_{1}N$ does not commute with $e_{2}N$. However, $EN/N$ is abelian. This is also a contradiction. Similarly, for the third case, one can use  (\ref{eq4.1})--(\ref{eq4.4}) to show that $e_{1}N$ does not commute with $e_{2}N$, also deriving a contradiction. Hence, the claim follows and  $EN/N$ is a cyclic group. 
The proof is completed.
\end{proof}

\section*{Acknowledgments} 

This research was partially supported by NNSF of China (No. 11401087).
The authors would like to thank Urban Jezernik and Primoz Moravec  for their help.

%%%%%%%%%%%%%%%%%%%%%%%%%%%%%%%%%%%%%%%%%%%%%%%%REFERENCES

%%%%%%%%%%%%%%%%%%%%%%%%%%%%%%%%%%%%%%%%%%%%%%%%%%%THE END
\end{document}